\documentclass{amsart}
\usepackage{latexsym}
\usepackage{amssymb}

\newcommand{\A}{$\mathfrak A$\;}
\newcommand{\Am}{$A$-$\mathfrak A$-module\;}
\newcommand{\AtB}{$A\hat{\otimes}_{\mathfrak A} B$\;}
\newcommand{\AtA}{$A\hat{\otimes}_{\mathfrak A} A$\;}
\newcommand{\D}{$D:A\rightarrow X^*$\;}
\newcommand{\am}{amenable\;}
\newcommand{\ma}{module amenable\;}
\newcommand{\md}{module derivative\;}
\newcommand{\mvd}{module virtual diagonal\;}
\newcommand{\mad}{module approximate diagonal\;}
\newcommand{\B}{Banach algebra\;}
\newcommand{\ba}{Banach $\mathfrak A$-module\;}

\newcommand{\bai}{bounded approximate identity\;}

\newcommand{\al}{\alpha . (a.x)}
\newcommand{\all}{(\alpha .a).x}
\newcommand{\alf}{\alpha . (a.f)}
\newcommand{\allf}{(\alpha .a).f}

\newcommand{\alffx}{(\alpha .f)(x)=f(x.\alpha)}
\newcommand{\affx}{(a.f)(x)=f(x.a)}
\newcommand{\alffb}{(\alpha .f)(b)=f(b.\alpha)}
\newcommand{\affb}{(a .f)(b)=f(ba)}
\newcommand{\alFFf}{(\alpha .F)(f)=F(f.\alpha)}
\newcommand{\aFFf}{(a .F)(f)=F(f.a)}

\newcommand{\aax}{\quad(\alpha\in\mathfrak A , a\in A, x\in X)}
\newcommand{\afx}{\quad(a\in A , f\in X^*, x\in X)}

\newcommand{\aalfx}{\quad(\alpha\in\mathfrak A , a\in A, f\in X^*, x\in X)}

\newcommand{\aalfb}{\quad(a,b\in A, \alpha\in\mathfrak A, f\in A^*)}
\newcommand{\aalfF}{\quad(a\in A, \alpha\in\mathfrak A , f\in A^*, F\in A^{**})}

\newcommand{\ts}{topological space\;}

\newtheorem{defi}{Definition}[section]
\newtheorem{prop}{Proposition}[section]
\newtheorem{theo}{Theorem}[section]
\newtheorem{lemm}{Lemma}[section]
\newtheorem{cor}{Corollary}[section]
\newtheorem{ex}{Example}[section]

\begin{document}

\title[Module Ameanability]{Module Amenability for semigroup algebras}
\author[M. Amini]{Massoud Amini}
\address{Department of Mathematics and Statistics\\ University of Saskatchewan
\\106 Wiggins Road, Saskatoon\\ Saskatchewan, Canada S7N 5E6\\mamini@math.usask.ca}
\thanks{Partially supported by a grant from University of Saskatchewan}
\maketitle

\begin{abstract}

We extend the concept of amenability of a Banach algebra $A$ to the case that there is an extra \A -module structure on $A$, and show that when $S$ is an inverse semigroup with subsemigroup $E$ of idempotents, then $A=\ell^1(S)$ as a Banach module over \A$=\ell^1(E)$ is module amenable iff $S$ is amenable. When $S$ is a discrete group, $\ell^1(E)=\mathbb C$ and this is just the celebrated Johnson's theorem. 

\end{abstract}

\section{Introduction}

The celeberated Johnson's theorem (in the discrete case) asserts that a discrete group $G$ is amenable if and only if the Banach algebra $\ell^1(G)$ is amenable. This fails to be true for discrete semigroups (even for the good cases like Clifford semigroups). Indeed Dunford and Namioka have shown that for a wide class of inverse semigroups (the class of $E$-unitary inverse semigroups) $\ell^1(S)$ is not amenable if the subsemigroup $E=E_S$ of idempotent elements is infinite [DN]. 

The concept of Johnson's amenability for Banach algebra has been a main stream in the theory of Banach algebras in the last fifty years. For some classes of operator algebras, however there have been some parallel concepts, among which one could mention the concept of {\it central amenability} for $C^*$-algebras [L2],[PR]. Also recently some research has been done on the {\it relative amenability} of Banach algebars [L1],[L3].

Here we develope the concept of {\it module amenability} for a class of Banach algebras which could somehow be considered as a generalization of all of the above approaches. In particular we apply this idea to the above mentioned problem and show that if $\ell^1(S)$ is considered appropriately as a $\ell^1(E)$-module, then its module amenability is equivalent to amenability of $S$, restoring the Johnson's theorem for the case of inverse semigroups.

The next section is devoted to the general theory of module amenability for Banach algebras in which we prove the analogues of the classical results on amenability of Banach algebras. Our main reference is [P], which in most cases we adapt almost the same proof. The details of proofs are given for the sake of completeness. In the last section, these results are used to prove the above mentioned version of the Johnson's theorem for inverse semigroups. We believe that this theory could well be applied to the case of topological (or measure) groupoids to get similar amenability results.

\section{module amenability}
Let \A be a \B  and $A$ be a \ba such that it has an associative product which makes it a \B which is compatible with the module action in the sense that 
$$
\alpha . (ab) =(\alpha .a)b , \quad (\alpha\beta).a=\alpha .(\beta .a) \quad (\alpha,\beta\in\mathfrak A , a,b\in A)
$$
and the same for the right action.

\begin{defi}
$A$ is called \ma (as an \A -module, if ambiguity may happen) if for any Banach space $X$ which is at the same time a Banach $A$-module and a \ba with the compatibility of actions
$$
(a.x).\alpha=a.(x.\alpha),\quad \al =\all \aax
$$
and the same for the other side actions, and each bounded map \D with 
$$
D(a+b)=D(a)+D(b), D(ab)=a.D(b)+D(a).b \quad\quad (a,b\in A),
$$
and
$$
D(\alpha .a)=\alpha .D(a),\quad D(a. \alpha)=D(a). \alpha \quad\quad (\alpha\in \mathfrak A ,a\in A),
$$
there is an $x\in X^*$ such that $D(a)=a.x-x.a=:D_x(a)\quad(a\in A)$. Note that $D$ is not assumed to be $\mathbb C$-linear and so it is not necessarily an \A -module homomorphism.
\end{defi}

Note that $X^*$ is also a Banach module over $A$ and \A with compatible actions under the canonical actions of $A$ and \A
$$
\alf =\allf \afx ,
$$
and the same for the right actions. Here the canonical actions of $A$ and \A on $X^*$ are defined by
$$
\alffx  , \quad \affx \aalfx ,
$$
and the same for right actions. 

We call $A$-modules $X$ which have a compatible \A -action as above, $A$-\A -{\it modules}, and derivations like $D$ in the above definition, 
the {\it module derivations}. Hence the above assertion is to say that if $X$ is an $A$-\A -module, then so is $X^*$ under the canonical actions. 
Also we use the notation $Z_{\mathfrak A}(A,X^*)$ for the set of all module derivations \D , and $B_{\mathfrak A}(A,X^*)$ for those which are inner 
and $H^1_{\mathfrak A}(A,X^*)$ for the quotient group (which we call the first relative (to \A) cohomology group of $X^*$). Hence $A$ is \ma iff 
$H^1_{\mathfrak A}(A,X^*)={0}$, for each $A$-\A -module $X$.

\bigskip
{\bf Blanket Assumption:}{\it All over this paper we fix $A$ and \A as above , and use notations $X$ and $D$ for arbitrary module and derivation as 
in the above definition, unless they are otherwise specified explicitly.}
\begin{prop}
If \A has a \bai , then amenability of $A$ implies its module amenability.
\end{prop} 
{\bf Proof} By the Cohen's factorization theorem [DW] for \A -modules $A$ and $X^*$, for each $a\in A, x\in X^*$ there are $\beta ,\gamma\in\mathfrak A$, 
$b\in A$, and $y\in X^*$ such that $a=\beta .b$ and $x=\gamma .y$. Therefore if $\{\alpha_i\}$ is a \bai in \A , then
\begin{align*}
D(\lambda a)&=D(\lambda(\beta .b)) =\lim_{i} D(\lambda\alpha_i .a)
 =\lim_{i} \lambda\alpha_i .D(a)=\lim_{i} \lambda\alpha_i .(\gamma .y)\\
 &=\lambda(\gamma .y)=\lambda D(a)
\end{align*}
for each $\lambda\in\mathbb C , a\in A$. Hence $D$ is $\mathbb C$-linear, and so inner.\qed

\vspace{.3 cm}
As we will see later in section 3, the converse is false. Hence module amenability is somehow weaker than amenability. Indeed another example 
in section 3 shows that module amenability even does not imply the existence of a bounded approximate identity. However we have the following weaker 
notion which is implied by module amenability.

\begin{defi} A bounded net $\{a_i\}$ in $A$ is called a bounded approximate commutator if 
$$
\lim_{i}\| a_i a-aa_i\| =0 \quad (a\in A).
$$
\end{defi}

Clearly each bounded approximate identity is a bounded approximate commutator.

\begin{prop}
If $A$ is \ma then it has a bounded approximate commutator .
\end{prop}
{\bf Proof} Consider $X=A^*$. Then $X$ and $X^*$ are $A$-\A -modules under the canonical actions
$$
\affb , \alffb \aalfb ,
$$
and
$$
\aFFf , \alFFf  \aalfF ,
$$ 
with right actions defined similarly. Now consider the canonical embedding 
$\theta :A\hookrightarrow A^{**}$, then clearly $\theta\in Z_{\mathfrak A}(A,X^*)$, so there is $x\in A^{**}$ such that $\theta(a)=a.x-x.a\quad (a\in A)$. 
Take any norm bounded (by $C$) net $\{a_i\}$ in $A$ so that $\theta(a_i)\rightarrow x$ in $\sigma(A^{**},A^*)$, then it is easy to see that 
$a_i a-aa_i\rightarrow 0$ for each $a\in A$, in $\sigma(A,A^*)$. Hence given $\epsilon > 0$ and finite subset $F\subseteq A$, there is a convex combination 
$a_{F,\epsilon}$ of elements of $F$, norm bounded by $C$, such that $\|a_{F,\epsilon} a -aa_{F,\epsilon}\|<\epsilon$, for each $a\in F$. Then 
$\{a_{F,\epsilon}\}_{(F,\epsilon)}$ forms a bounded approximate commutator for $A$.  \qed

\vspace{.3 cm}
Given \A -modules $A$ and $B$, let \AtB be the projective module tensor product of $A$ and $B$ [R]. This is the quotient of the usual projective tensor product 
$A\hat{\otimes} A$ by the closed ideal $I$ generated by elements of the form $z.a\otimes b -a\otimes z.a$ for $z\in\mathfrak A, a,b\in A$. If $A$ and $B$ are also 
Banach algebras with compatible action, then so is \AtB . Also $(A\hat{\otimes}_{\mathfrak A} B)^*\cong {\mathfrak L}_{\mathfrak A}(B,A^*)$ , where the right hand side 
is the space of all \A -module morphisms from $B$ to $A^*$ [R]. In particular \AtA is a Banach $A$-\A -module. Here the second $A$ in the tensor product is understood to be 
with with the opposite product. Consider $\omega \in {\mathfrak L}(A\hat{\otimes} A, A)$ defined by $\omega(a\otimes b)=ab \quad (a,b\in A)$. Then both $\omega$ and its dual 
conjugate $\omega^{**} \in {\mathfrak L}((A\hat{\otimes} A)^{**}, A^{**})$ are \B homomorphisms. Now as $I$ is an ideal 
of $A\hat{\otimes} A$, so $\omega(I)$ is an ideal of $A$, 
and if $J$ is the closure of $\omega(I)$ we can define $\tilde{\omega} : A\hat{\otimes}_{\mathfrak A} A=A\hat{\otimes} A/I\rightarrow A/J$ by 
$$
 \tilde{\omega}(a\otimes b +I)=ab+J\quad(a,b\in A).
$$
This extends to an element $\tilde\omega\in {\mathfrak L}(A\hat{\otimes}_{\mathfrak A} A, A/J)$ and both $\tilde\omega$ and its dual conjugate 
$\tilde\omega^{**} \in {\mathfrak L}((A\hat{\otimes}_{\mathfrak A} A)^{**}, A^{**}/J^{\perp\perp})$ are $A$-\A-module homomorphisms. 

\begin{defi}
A bounded net $\{\tilde u_i\}$ in $A\hat{\otimes}_{\mathfrak A} A$ is called a module approximate diagonal if ${\tilde\omega(\tilde u_i}$ is a 
(bounded) approximate identity of $A/J$ and $\lim_i \|\tilde u_i .a-a.\tilde u_i\|=0\quad (a\in A)$. An element 
$\tilde M\in (A\hat{\otimes}_{\mathfrak A} A)^{**}$ is called a module virtual diagonal if 
$$
\tilde\omega^{**}(\tilde M).a=a, \quad \tilde M.a=a.\tilde M\quad (a\in A).
$$
\end{defi}

\begin{prop} The following are equivalent:
 
$(i)$\,\, $A$ has a \mvd

$(ii)$ There is $M\in (A\hat{\otimes} A)^{**}$ such that 
$$
\omega^{**}(M).a-a\in J^{\perp\perp}, \quad M.a-a.M\in I^{\perp\perp}\quad (a\in A),
$$
where $I$ is the closed ideal of $A\hat{\otimes} A$ generated by elements of the form $z.a\otimes b -a\otimes z.a$ for $z\in\mathfrak A, a,b\in A$, 
and $J=\bar{\omega(I)}$.

In particular if $A$ has a virtual diagonal, then it has a \mvd .
\end{prop}
{\bf Proof} If $M$ is as in $(ii)$, define $\tilde M\in (A\hat{\otimes}_{\mathfrak A} A)^{**}=(A\hat{\otimes} A)^{**}/I^{\perp\perp}$ 
by $\tilde M=M+I^{\perp\perp}$. Then given $a\in A$, since $M.a-a.M\in I^{\perp\perp}$
, clearly $\tilde M.a=a.\tilde M$, and since $\omega^{**}(M).a-a\in J^{\perp\perp}$ and 
$\tilde\omega^{**}(\tilde M)=\tilde\omega^{**}(M+ I^{\perp\perp})=\omega^{**}(M)+J^{\perp\perp}$, clearly $\tilde\omega^{**}(\tilde M).a=a$.

Conversely, if $\tilde M$ is a \mvd , then choose any $M\in (A\hat{\otimes} A)^{**}$ such that $\tilde M=M+I^{\perp\perp}$, 
then given $a\in A$, $ M.a-a.M + I^{\perp\perp}=\tilde M.a-a.\tilde M=0\in (A\hat{\otimes} A)^{**}/I^{\perp\perp}$, so  
$ M.a-a.M \in I^{\perp\perp}$. Also $\omega^{**}(M).a-a +J^{\perp\perp}=\tilde\omega^{**}(M).a-a=0\in A^{**}/J^{\perp\perp}$, 
so  $\omega^{**}(M).a-a \in J^{\perp\perp}$.\qed

\vspace{.3 cm}
The following is proved similarly.

\begin{prop} The following are equivalent:
 
$(i)$\,\, $A$ has a \mad

$(ii)$ There is a bounded net $\{u_i\}\in (A\hat{\otimes} A)$ such that for each $a\in A$ the nets  
$\{\omega(u_i).a-a\}$ and $\{u_i .a -a.u_i\}$ converge to an element of $J$ and $I$, respectively.

In particular if $A$ has an approximate diagonal, then it has a \mad . \qed
\end{prop}

\begin{lemm} With above notations

$(i)$\,\, If $A$ is unital with unit $e$, then $H^1_{\mathfrak A}(A,X)\cong H^1_{\mathfrak A}(A,eXe)$,

$(ii)$\,  $H^1_{\mathfrak A}(A\oplus \mathfrak A,X)\cong H^1_{\mathfrak A}(A,X)$,

$(iii)$ If $A$ has \bai , then $H^1_{\mathfrak A}(A,X^*)\cong H^1_{\mathfrak A}(A,(AXA)^*)$,

$(iv)$ If \A has \bai , then $H^1_{\mathfrak A}(A,X^*)\cong H^1_{\mathfrak A}(A,(\mathfrak A X\mathfrak A)^*)$.
\end{lemm}
{\bf Proof} $(i)$ Consider the $A$-\A -module homomorphisms $id ,\ell,r:X\rightarrow X$ defined by $id(x)=x$, $\ell(x)=e.x$ and $r(x)=x.e$. Then 
$$
X\cong eXe \oplus (id-r)\circ\ell (X)  \oplus (id-\ell)\circ r(X)
 \oplus (id-r)\circ(id-\ell)(X),
$$
with pairwise commuting projections. Now it is easy to check that each summand of the right hand side is a Banach $A$-\A -module and that the 
first relative cohomology groups of the last three terms vanish.

$(ii)$ Clearly $X$ is an $(A\oplus\mathfrak A)$-\A -module. Let $D$ be any module derivation from $A\oplus \mathfrak A$ to $X$, then under the 
identification \A$\cong\{0\}\oplus$\A , for each $a\in A$ and $\alpha ,\beta\in\mathfrak A$ we have
$$
(a,\beta).D(0,\alpha)=D((a,\beta).(0,\alpha))-D((a,\beta)).\alpha=D((a,\beta).\alpha)-D((a,\beta)).\alpha =0.
$$
Hence $D\circ\theta$ is a module derivation from $A$ to $X$, for the canonical embedding $\theta:A\rightarrow A\oplus\mathfrak A$. It is easy 
to see that $D\mapsto D\circ \theta$ lifts to an isomorphism of the corresponding relative cohomology groups.

$(iii)$ Let $X_1 =AXA$ and $X_2 =AX$. As $A$ contains a \bai $\{a_i\}$, Cohen's factorization theorem shows that $X_1$,{$X_2$ are $A$-\A -modules. 
Now $A$ has a zero left action on $X_2^\perp \cong (X/X_2)^*$, 
so $H^1_{\mathfrak A}(A,X^\perp _2)={0}$. Also $\{a_i\}$ is a 
bounded net in $\mathfrak L_{\mathfrak A}(X^*)\cong(X\hat{\otimes}_{\mathfrak A} X^*)^*$, so 
passing to a subnet we may assume that it is $w^*$-convergent to some $F\in\mathfrak L(X^*)$. Then 
$$
<a.x,F(f)>=\lim_i<a.x,f.a_i>=\lim_i <a_i a.x,f>=<a.x,f> ,
$$
for each $a\in A, f\in X^*, x\in X$.
Hence $I-F$ is a projection of $X^*$ onto $X^\perp _2$, where $I$ is the identity map on $X^*$. 
Therefore $X^*\cong X_2\oplus X^\perp _2$ and so $H^1_{\mathfrak A}(A,X^*)\cong H^1_{\mathfrak A}(A,X^*_2)$. 
Now $H^1_{\mathfrak A}(A,X^*_1)\cong H^1_{\mathfrak A}(A,X^*_2)$ is similar.

$(iv)$ Proof is quite similar to $(iii)$. \qed

\begin{theo}
The following are equivalent

$(i)$\;\; $A$ is \ma ,

$(ii)$\; $A$ has a \mad ,

$(iii)$ $A$ has a \mvd .
\end{theo} 
{\bf Proof} The equivalence of $(ii)$ and $(iii)$ follows exactly like in the classical case [P, Lemma 1.6]. 
Assume that $A$ has a \mvd $\tilde M\in (A\hat{\otimes}_{\mathfrak A} A)^{**}$ and let $M\in (A\hat{\otimes} A)^{**}$ be 
the corresponding element as in $(ii)$ of Proposition 2.3. We may assume that $\tilde M$ is a $w^*$-limit point of a \mad $\{\tilde u_i\}$. 
Take any Banach $A$-\A -module $X$ and any module derivation \D . By above lemma, we may assume that $X$ is essential, 
that is $X=AXA$. To each $x\in X$, there corresponds $F_x\in (A\hat{\otimes} A)^*$ via
$$
F_x(a\otimes b)=<x, a.D(b)> \quad (a,b\in A).
$$
For each $a\in A, x\in X$. Then $D(a).x$ could thought as an element of $A^*$ via 
$$
<D(a).x,b>=<D(a),b.x>\quad (b\in A),
$$
and it is easy to see that 
$$
F_{(a.x-x.a)}=a.F_x -F_x .a +\omega^*(D(a).x)\quad (a\in A, x\in X),
$$
(see proof of [P,1.7]). Putting $\tilde F_x=F_x +I^\perp$ and $f(x)=<F_x,M>\quad(x\in X)$, we have 
$\tilde F_x\in (A\hat{\otimes}_{\mathfrak A} A)^{*}$ and $f\in X^*$, it is easy to check that 
$$
\tilde F_{(a.x-x.a)}=a.\tilde F_x -\tilde F_x .a +\tilde\omega^*(D(a).x+J^\perp)\quad (a\in A, x\in X),
$$
and so
\begin{align*}
<x,D_f(a)> & =<\tilde F_{(a.x-x.a)},\tilde M>\\
&=<a.\tilde F_x -\tilde F_x .a ,\tilde M>+<\tilde\omega^*(D(a).x+J^\perp),\tilde M>\\
& =<\tilde F_x,\tilde M .a-a.\tilde M> +\lim_i<x.\tilde\omega(\tilde u_i),D(a)>=<x,D(a)>,
\end{align*}
where the last equality is because $X$ is essential. Hence $D=D_f$. 

Conversely assume that $A$ is \ma ; then by Proposition 2.2 , it has a \bai $\{a_i\}$. Put $\dot{a}_i=a_i +J\in A/J$, then passing 
to a subnet, we may assume that $\{\dot a_i\otimes \dot a_i\}$ is $w^*$-convergent to an element $N\in (A\hat{\otimes}_{\mathfrak A} A)^{**}$.  
Clearly $\tilde\omega$ 
vanishes on the range of $D_N$ and so $D_N$ could be regarded as a 
module derivation of $A$ into the $A$-\A -module $K=ker(\tilde\omega^{**})$. 
By module amenability of $A$, there is $N^{'}\in K$ such that $D_N=D_{N^{'}}$, so $M=N-N^{'}$ is a \mvd .\qed

\begin{prop} If $A$ and $B$ are Banach algebras and Banach \A -modules with compatible actions, and there is a continuous Banach algebra homomorphism and 
module morphism $\varphi$ from $A$ onto a dense subset of $B$, and $A$ is \ma , then so is $B$.
\end{prop}
{\bf Proof} If $X$ is any $B$-\A -module then via $\varphi$, $X$ could be regarded as an $A$-\A -module, and each module derivations $D:B\rightarrow X^*$, 
gives a module deivation $D\circ\varphi :A\rightarrow X^*$, which is inner. By density of 
range of $\varphi$ and continuity of $D$, then $D$ would be inner.\qed

\vspace{.3 cm}
Next assume that $A$ has a \bai ${a_i}$ and consider the algebra of \A -multipliers of $A$
$$
M_{\mathfrak A}(A)=\{(T_1,T_2): T_1, T_2\in \mathfrak L_{\mathfrak A}(A): T_1(ab)=T_1(a)b, T_2(ab)=aT_2(b) \;(a,b\in A)\}.
$$
Then $M_{\mathfrak A}(A)$ is an $A$-\A -module and $A$ embeds in $M_{\mathfrak A}(A)$ via $a\mapsto (S_a,T_a)$, where 
$S_a(b)=ab, T_a(b)=ba\quad(a,b\in A)$. For any element $T=(T_1,T_2)$ 
of $M_{\mathfrak A}(A)$ it is easy to see that 
$\|T_1\|=\|T_2\|$ and if we put $\|T\|$ equal to this common value, then $M_{\mathfrak A}(A)$ becomes a Banach \Am . Also for each \Am $X$, 
$M_{\mathfrak A}(A)$ acts on $X$ via 
$$
T.x=\lim_i T_1(a_i).x ,\quad x.T=\lim_i x.T_2(a_i)\quad\quad(x\in X, T=(T_1,T_2)\in M_{\mathfrak A}(A)),
$$
which makes $X$ a $M_{\mathfrak A}(A)$-\A -module. Also given a module derivation $D:M_{\mathfrak A}(A)\rightarrow X^*$, the restriction $D^{'}$ of 
$D$ to $A$ is a module derivation on $A$. 

\begin{prop} 
With the above notation, if $A$ has a \bai , then $H^1_{\mathfrak A}(M_{\mathfrak A}(A),X^*)\cong H^1_{\mathfrak A}(A,X^*)$.
\end{prop}
{\bf Proof} It is enough to show that if the right hand side is $\{0\}$, then so is the left hand side. If $D$ and $D^{'}$ are as above, then $D^{'}=D_x$, 
for some $x\in X^*$. Put $D_i=T_1(a_i).x-x.T_2(a_i)$, 
for each index $i$ and $T=(T_1,T_2)\in  M_{\mathfrak A}(A)$, then the net ${D_i}$ is 
uniformely norm bounded, and so, passing to a subnet, we may assume that it converges to some $D_0\in{\mathfrak L}_{\mathfrak A} (M_{\mathfrak A}(A), X^*)$ 
in $\sigma(X^*, M_{\mathfrak A}(A))$. Then it is easy to show that $D_0=D$ and since clearly $D_0=D_x$ on $A$, so by the continuity of these derivations in 
the strict topology of $M_{\mathfrak A}(A)$ 
and density of $A$ in $M_{\mathfrak A}(A)$ in the strict topology, we get $D=D_x$ on $M_{\mathfrak A}(A)$.\qed

\begin{cor} Let $J$ be a closed ideal of $A$ which is \A -invariant, i.e. $\mathfrak A  .J\subseteq J$. If $J$ has a \bai and $A$ is \ma then $J$ is 
also \ma .
\end{cor}
{\bf Proof} Let $X$ be an essential Banach $J$-\A -module and $D:M_{\mathfrak A}(J)\rightarrow X^*$ be a bounded module derivative. By definition of 
$M_{\mathfrak A}(J)$, 
there is an \A -module morphism $\varphi:A\rightarrow M_{\mathfrak A}(J)$ 
and $D\circ\varphi$ is a module derivation on $A$, so it is inner. Hence $D$ is inner on $J$. By 
the same argument as above $D$ is inner on $M_{\mathfrak A}(J)$ and so the result follows from the above proposition.\qed

\begin{cor}
If $J$ is an \A -invariant closed ideal of $A$ with \bai , then $A$ is \ma iff both $J$ and $A/J$ are \ma .
\end{cor}
{\bf Proof} One direction is trivial, for the other, assume that both $J$ and $A/J$ are \ma and let \D be any bounded module derivative, then the 
restriction $D^{'}$ of $D$ to $J$ is inner. 
It is easy to see that the range of $D-D^{'}$ 
vanishes on both $JX$ and $XJ$. Now The close submodule $X_J$ of $X$ generated by $JX\cup XJ$ is a 
Banach \Am and $D-D^{'}:A\rightarrow \hat{X}_J\cong(X/X_J)^*$ is a module derivative which vanishes on $J$, 
so it could be regarded as a module deivative on $A/J$, and so it is inner. Hence $D$ is inner.\qed

\begin{defi}
Given a net of Banach algebras and \A-modules $\{A_i\}$ with compatible actions, 
we say that they are simultaneously module amenable if there is a constant $c>0$ such that for each index $i$, each $A_i$-$\mathfrak A$-module 
$X_i$, each module derivation $D_i:A_i\rightarrow X^*_i$, 
there is $x_i\in X^*_i$ such that $D_i=D_{x_i}$ and $\|x_i\|\leq c\|D_i\|$.
\end{defi}

\begin{prop} If a direct system $\{A_i\}$ is simulteneously module amenable, then $A:=\varinjlim A_i$ is \ma .
\end{prop}
{\bf Proof} We may assume that $A_i\subseteq A_i\subseteq A\quad(i\leq i)$ and $A=(\cup_{i} A_i)^{\bar{}}$. Let $X$ be a Banach $A$-\A -module 
and \D is a bounded \md ,then the restriction $D_i$ of $D$ to $A_i$ is a \md on $A_i$, so there is $c>0$ and $z_i\in X^*$ such that $D=D_{z_i}$ 
on $A_i$ and $\|z_i\|\leq c\|D\|$. Passing to a 
subnet we may assume that $\{z_i\}$ is $w^*$-convergent to some $z\in X^*$. Take any $a\in \cup_i A_i$,then $a\in A_i\quad(i>i_{0})$, for some 
index $i_0$. Then
\begin{align*}
<D_z(a),x>&=<z,x.a-a.x>=\lim_i <x.a-a.x,z_i>\\
&=\lim_i<D_{z_i}(a),x>=<D(a),x>\quad(a\in A, x\in X^*),
\end{align*}
so $D=D_z$ on a dense subset of $A$, and so on $A$. \qed

\vspace{.3 cm}
Next we consider the module amenability of the module tensor product.


\begin{prop}
If $A$ and $B$ are module amenable then so is \AtB .
\end{prop}
{\bf proof} As $A$ and $B$ admit \bai by Proposition 2.2, so does \AtB . Take any essential Banach \AtB-\A-module $X$. The mappings $a\mapsto \sigma_a$ 
and $b\mapsto\tau_b$ defined by 
$$
\sigma_a(c\otimes d)=ac\otimes d ,\quad \tau_b(c\otimes d)=c\otimes bd\quad (a,c\in A, b,d\in B),
$$
extend and then lift to commuting continuous Banach algebra homomorphisms and \A-module morphisms of $A$ and $B$ onto closed subalgebras $A_1$ and $B_1$ 
of $M_{\mathfrak A}(A\hat{\otimes}_{\mathfrak A} B)$, 
and $X$ is then a Banach $A_1$-\A-module and $B_1$-\A-module. Now each \md $D:A\hat{\otimes}_{\mathfrak A} B \rightarrow X^*$ gives rise to 
some \md $D^{'}:M_{\mathfrak A}(A\hat{\otimes}_{\mathfrak A} B)\rightarrow X^*$, whose restriction $D_1$ to $B_1$ is inner, say $D^{'}=D_x$ on $B_1$, 
for some $x\in X^*$. Consider $d:X\rightarrow{\mathfrak L}_{\mathfrak A} (B_1,X)$, which sends each $y\in X$ to the restriction of $D_y$ to $B_1$. 
Then the fact that $D^{'}-D_x$ is zero on $B_1$ implies that $D^{'}-D_x$ is sending $A_1$ into $Im(d)^\perp \cong(X/Im(d))^*$ (since $\sigma$ 
and $\tau$ commute). Now $Im(d)$ and so $Y:=X/Im(d)$ is an $A_1$-\A-module (again since  $\sigma$ and $\tau$ commute), and so $D^{'}-D_x=D_y$ 
on $A_1$, for some $y\in Y$. But both sides of the last equality are equal zero on $B_1$, and so by the fact that $A_1\cup B_1$ is dense 
in $M_{\mathfrak A}(A\hat{\otimes}_{\mathfrak A} B)$, $D^{'}$ is inner on $M_{\mathfrak A}(A\hat{\otimes}_{\mathfrak A} B)$, so $D$ is 
inner on $A\hat{\otimes}_{\mathfrak A} B$.\qed

\begin{prop}
If \A is unital and $B$ is a unital \am \B , then $A:=B\hat{\otimes}_{\mathfrak A} \mathfrak A $ is \ma . In particular, if \A is unital then it is \ma 
as a module over itself.
\end{prop}
{\bf Proof} As \A and $B$ are unital, $B$ and \A could be identified with a subalgebra of $A$. Given $A$-\A-module $X$, and bounded \md \D , $X$ 
is also a $B$-\A-module and the restriction $D_1$ of $D$ to $B$ is a \md . As $B$ is unital, $D$ is a derivation, and so inner. Also clearly the 
restriction $D_2$ of $D$ to \A is the zero map (as \A is unital). Now for each $b\in B, \alpha\in\mathfrak A$,
\begin{align*}
D(b\otimes\alpha)&=(1\otimes\alpha)D_1(b)+D_2(\alpha).b=(1\otimes\alpha)(b.x-x.b)+0\\
&=(b\otimes\alpha).x-x.(b\otimes\alpha)=D_x(b\otimes\alpha) .
\end{align*}
Now as \A is unital, $D$ is indeed linear, hence $D=D_x$ on $B\hat{\otimes} \mathfrak A$.\qed

\vspace{.3 cm}
Using this result we can now provide a family of examples of module amenable Banach algebras, Here are some:

\begin{ex}
$(i)$ \;\;If $\Omega$ is a compact \ts , then  $C(\Omega, \mathfrak A)\cong C(\Omega)\hat{\otimes} \mathfrak A$ is \ma .

$(ii)$ \;If $\mathfrak H$ is a separable Hilbert space, then 
$(\mathfrak K(\mathfrak H)\oplus\mathbb C I)\hat{\otimes}\mathfrak A\cong(\mathfrak K(\mathfrak H)\hat{\otimes}\mathfrak A)\oplus\mathfrak A$ is \ma .

$(iii)$ If $G$ is a discrete \am group, then $\ell^1(G)\hat{\otimes}\mathfrak A$ is \ma , in particular, for $\mathfrak A=\ell^1(H)$, 
where $H$ is a subgroup of $G$, $\ell^1(G\times H)$ is \ma as an $\ell^1(H)$-module.
\end{ex}

\section{module amenability of semigroup algebras}

In this section we consider an important example which was the motivation for writing this paper. We consider an inverse semigroup $S$ with 
idempotents $E$ and show that $\ell^1(S)$ is \ma as an $\ell^1(E)$-module if and only if $S$ is amenable.

\begin{defi} A discrete semigroup $S$ is called an inverse semigroup if for each $x\in S$ there is a unique element $x^*\in S$ such that 
$xx^*x=x$ and $x^*xx^*=x^*$. An element $e\in S$ is called an idempotent if $e=e^*=e^2$. The set of idempotent elements of $S$ is denoted by $E$.
\end{defi}

\begin{defi} A discrete semigroup $S$ is called amenable if there is an invariant mean on $\ell^\infty(S)$, namely an element $m\in\ell^\infty(S)^*$ such that $m(1)=\|m\|=1$ and $m(s.f)=m(f.s)=m(f)\quad(s\in S,f\in \ell^\infty(S))$, where 
$$
f.s(t)=f(ts),\quad s.f(t)=f(st)\quad(s,t\in S,f\in \ell^\infty(S)).
$$
\end{defi}

\vspace{.3 cm}
It is easy to see that $E$ is indeed a commutative subsemigroup of $S$. In particular $\ell^1(E)$ could be regarded as a subalgebra of $\ell^1(S)$, and thereby $\ell^1(S)$ is a Banach algebra and a Banach $\ell^1(E)$ module with compatible action. The same then would be true for $\ell^\infty(S)$. Of course one may change the action of $\ell^1(E)$ on $\ell^1(S)$ to get different module amenability results. One lesson we learned from the proof of Proposition 2.2 was that sometimes it is helpful to consider the action from one side to be some sort of trivial action (zero action in that case). We adapt that idea here and let $\ell^1(E)$ act from right on $\ell^1(S)$ by multiplication and as identity from left, that is
$$
\delta_e .\delta_s =\delta_s, \; \delta_s .\delta_e =\delta_se = \delta_s *\delta_e \quad (s\in S, e\in E).
$$

\begin{lemm}
With the above notation

$(i) \;\ell^1(S)\hat{\otimes}_{\ell^1(E)}\ell^1(S)\cong\ell^1(S\times S)/I$,
where $I$ is the closed ideal of $\ell^1(S\times S)$ which is closed linear span of the set of 
elements of the form $\delta_{(set,x)}-\delta_{(st,x)}$, where $s,t,x\in S$ and $e\in E$.

$(ii) (\ell^1(S)\hat{\otimes}_{\ell^1(E)}\ell^1(S))^{*}\cong \ell^\infty(S\times S)/I^\perp$,
where 
$$
I^\perp =\{f\in  \ell^\infty(S\times S): f(set,x)=f(st,x)\quad(s,t,x\in S, e\in E)\}.
$$
\end{lemm}
{\bf Proof} $(i)$ follows directly from the definition of the module projective tensor product. For $(ii)$, if $f\in I^\perp$, 
then for each $s,t,x\in S$ and $e\in E$,
\begin{align*}
0&=<f,\delta_{(set,x)}-\delta_{(st,x)}>=\sum_{u,v}f(u,v)(\delta_{(set,x)}(u,v)-\delta_{(st,x)}(u,v)) \\
& =f(set,x)-f(st,x).
\end{align*}
Conversely, if $f$ satisfies the given relation, then clearly $<f,u>=0$ for each $u\in I$ which is a finite linear combination of elements of the form  
$\delta_{(set,x)}-\delta_{(st,x)}$, where $s,t,x\in S$ and $e\in E$. By continuity, then $f\in I^\perp$.\qed

\vspace{.3 cm}
Consider $\omega : \ell^1(S)\hat{\otimes} \ell^1(S)=\ell^1(S\times S)\rightarrow\ell^1(S)$ defined by 
$$
\omega(g)(s,t)=g(st) \quad (s,t\in S, g\in\ell^1(S)),
$$
then
$$
\omega^*(h)(s,t)=h(st) \quad (s,t\in S, h\in\ell^\infty(S)).
$$
Now if 
$$
f.s(t,t^{'})=f(st,t^{'})\quad s.f(t,t^{'})=f(t,t^{'}s)\quad (s,t,t^{'}\in S, f\in\ell^\infty(S\times S)), 
$$
and for $M\in \ell^\infty(S\times S)^{**}$ and $f\in \ell^\infty(S\times S)$,
$$
M.s(f)=M(s.f), \quad s.M(f)=M(f.s),
$$
then $M$ is a virtual diagonal of $\ell^1(S)$ if and only if, for each $s\in S$
$$
M.s=s.M , \quad \omega^{**}(M).s=s,
$$
on $\ell^\infty(S\times S)$, where in the second equality the left hand side is defined by
$$
\omega^{**}(M).s(h)=M(s.(\omega^*(h)))\quad (s\in S, h\in \ell^\infty(S)),
$$
and the right hand side is the functional of evaluation at $s$ on $\ell^\infty(S))$ [DN]. Now if $I$ is as in the above lemma and $J=\omega(I)^{\bar{}}$, 
then by Proposition 2.3 , 
$M\in \ell^1(S\times S)^{**}=\ell^\infty(S\times S)^{*}$ gives rise to a \mvd for  $\ell^1(S)$ if and only if, for each $s\in S$ the equalities
$$
M.s=s.M , \quad \omega^{**}(M).s=s,
$$
hold on $I^\perp\subseteq\ell^\infty(S\times S)$ and $J^\perp=\omega(I)^\perp\subseteq \ell^\infty(S)$, respectively, where
$$
J^\perp=\omega(I)^\perp=\{h\in\ell^\infty(S): h(set)=h(st)\quad (s,t\in S, e\in E)\}.
$$
Such an element $M$ exists if and only if $\ell^1(S)$ is \ma .

\vspace{.3 cm}
Next consider the congruence relation $\backsim$ on $S$ defined by $s\backsim t$ if and only if there is $e\in E$ such that $se=te$. The quotient semigroup 
$G_S:=S/{\backsim}$ is then a group. It is indeed the 
maximal group homomorphic image of $S$ [Mu]. Also the inverse semigroup $S$ is amenable if and only if the discrete group $G_S$ is amenable [DN]. 

\begin{lemm}
With above notation, $\ell^1(G_S)$ is a quotient of $\ell^1(S)$ and so the above action of $\ell^1(E)$ on $\ell^1(S)$ 
lifts to an action of $\ell^1(E)$ on $\ell^1(G_S)$, 
making it a Banach  $\ell^1(E)$-module. 
\end{lemm}
{\bf Proof} Consider the quotient map $\pi :S\rightarrow G_S$, then extending  $\pi$ by linearity and noting that for each $n\geq 1$ , 
each $c_1,\dots ,c_n \in \mathbb C$, and each $s_1,\dots ,s_n\in S$
$$
\|\sum_{i=1}^{n} c_i \delta_{\pi(s_i)}\|_1 = \|\sum_{i=1}^{n} c_i \delta_{s_i}\|_1
$$
one can extend $\pi$ to a continuous Banach algebra epimomorphism
$\pi :\ell^1(S)\rightarrow \ell^1(G_S)$. \qed

\begin{lemm}
With the setup of above lemma, $\ell^1(G_S)$ is \ma if and only if it is amenable.
\end{lemm}
{\bf Proof} The left action of $\ell^1(E)$ on $\ell^1(S)$ and so on $\ell^1(G_S)$ is trivial. As for the right action, we have 
$$
\delta_{\pi(s)} .\delta_e =\delta_{\pi(se)}=\delta_{\pi(s)} \quad(s\in S, e\in E).
$$
Hence the right action is also trivial and so 
$$
\ell^1(G_S)\hat{\otimes}_{\ell^1(E)}\ell^1(G_S)\;\cong\;\ell^1(G_S)\hat{\otimes}\ell^1(G_S),
$$
and the result follows from Proposition 2.3.\qed

\vspace{.3 cm}
Now we are ready to prove the main result of this section.

\begin{theo}
Let $S$ be an inverse semigroup with idempotents $E$. Consider $\ell^1(S)$ as a Banach module over $\ell^1(E)$ with the multiplication 
right action and the trivial left action. Then $\ell^1(S)$ is \ma if and only if $S$ is amenable.
\end{theo}
{\bf proof} If $\ell^1(S)$ is \ma then so is $\ell^1(G_S)$ by Lemma 3.2 and Proposition 2.5. Hence $\ell^1(G_S)$ is amenable by above Lemma, 
and so $G_S$ is amenable by Johnson's theorem. Therefore $S$ is amenable by [DN, thm.1].

Conversely if $\mu$ is a right invariant mean on $S$ and $M$ is defined on $\ell^\infty(S\times S)$ by
$$
M(f)=\int_{S} f(s^*,s)d\mu(s) ,
$$
then $M$ is clearly a bounded linear functional and $M(1\otimes 1)=\mu(1)=1$. Also for each $s\in S$ and $f\in \ell^\infty(S\times S)$,
\begin{align*}
s.M(f)=M(f.s) &=\int_{S} f(st^*,t)d\mu(t)\\
 &=\int_{S} f(s(ts)^*,ts)d\mu(t)\\
 &=\int_{S} f(ss^*t^*,ts)d\mu(t)\\
 &=\int_{S} f((tss^*)^*,(tss^*)s)d\mu(t)\\
 &=\int_{S} f(t^*,ts)d\mu(t)\\
 &=M(s.f)=M.s(f),
\end{align*}
and for each $s\in S$ and $f\in J\subseteq\ell^\infty(S\times S)$,

\begin{align*}
\omega^{**}(M).s(f) &=\omega^{**}M(f.s) \\
 &=M(\omega^*(f.s))=\int_{S} \omega^*(f.s)(t^*,t)d\mu(t)\\
 &=\int_{S} f.s(t^*t)d\mu(t)=\int_{S} f(st^*t)d\mu(t)\\
 &= f(s)\int_{S} d\mu(t)=f(s),
\end{align*}
where the last equality between integrals is because by the discussion after Lemma 3.1, $f(se)=f(s)$ , for each $e\in E$.
Hence $M$ gives rise to a \mvd for $\ell^1(S)$ and so $\ell^1(S)$ is \ma . \qed

\vspace{.3 cm}
Now we are ready to give the two countraexamples mentioned in section 2. For the second example we use the fact that if $S$ satisfies the $D_k$ condition of Duncan and Namioka, 
for some positive integer $k$ if and only if $\ell^1(E)$ has a bounded approximate identity if and only if $\ell^1(S)$ has a bounded approximate identity[DN].

\begin{ex} 
$(i)$ There is an inverse semigroup $S$ for which $\ell^1(S)$ is \ma but not amenable .

$(ii)$ There is an inverse semigroup $S$ for which $\ell^1(S)$ is \ma but has no \bai .
\end{ex}
{\bf Proof} For the first example, take any amenable inverse semigroup $S$ with infinite number of idempotents (there are a lot of them among Clifford semigroups). 
Then $\ell^1(S)$ is \ma by above theorem, but not amenable [DN]. For the second example take any {\it Brandt semigroup} $S$ of an amenable group over an infinite 
index set then clearly $S$ is amenable and so again $\ell^1(S)$ is \ma by above theorem, but it has no \bai [DN].\qed

\vspace{.3 cm}
Next we can consider the module amenability of the {\it reduced} semigroup $C^*$-algebra $C^*_r(S)$ and the semigroup von Neumann algebra
VN(S) generated by $C^*_r(S)$. We refer the reader 
to [Pa] for definitions. Here we only need to know that $C^*_r(S)$ is a homomorphic image of $\ell^1(S)$[Pa]. 
Therefore we may consider $C^*_r(S)$ as an $\ell^1(E)$-module (with the induced actions from the action of $\ell^1(E)$ on $\ell^1(S)$ as in the 
above theorem). Also VN(S) is homomorphic image of the second conjugate of $C^*_r(S)$ and so carries induced action of $\ell^1(E)$ similarly. Now in this 
setting we have the following two partial results. 
\begin{cor}
With the setting of the above theorem, if $S$ is amenable then 
$C^*_r(S)$ is \ma .
\end{cor}
{\bf Proof} If If $S$ is amenable, then $C^*_r(S)$ is \ma by above theorem and Proposition 2.3.\qed

\begin{cor}
With the setting of the above theorem, if every maximal subgroup of $S$ is amenable and $S$ satisfies the $D_k$ condition of Duncan and Namioka, 
for some positive integer $k$, then $VN(S)$ is \ma .
\end{cor}
{\bf Proof} Under the first condition $VN(S)$ is amenable [Pa] and under the second condition $\ell^1(E)$ has a bounded approximate identity [DN]. Hence the result
follows from Proposition 2.1.\qed
  
\vspace{.5 cm}
{\small {\bf Acknowledgement}: This paper is prepared while the author was visiting the  department of Mathematics and Statistics of University of 
Saskatchewan. I would like to thank the hospitality I have received there. 
In particular I would like to thank Dr. M. Khoshkam and Dr. David Cowan for many stimulating discussions.

\end{document}